\documentclass[11pt,a4paper]{amsart}
\usepackage{a4wide}
\usepackage[utf8]{inputenc}
\usepackage[T1]{fontenc}
\usepackage{a4wide}
\newcommand{\Q}{\mathbb{Q}}
\newcommand{\2}{\mathbb{D}}

\newcommand{\Prop}{\mathsf{Prop}}
\newcommand{\Form}{\mathsf{Form}}
\newcommand{\model}[1]{\mathcal{#1}}
\renewcommand{\implies}{\rightarrow}
\newcommand{\struc}[1]{\langle #1 \rangle}
\newcommand{\Val}{\mathrm{Val}}
\newcommand{\lucas}{\mbox{\L}}
\newcommand{\logic}[1]{\mathsf{\mathbf{#1}}}
\newcommand{\K}{\mathsf{K}}
\newcommand{\proves}{\vdash}
\newcommand{\lang}{\mathcal{L}}
\newcommand{\var}[1]{\mathcal{#1}}
\newcommand{\classop}[1]{\mathbb{#1}}
\newcommand{\free}{\mathcal{F}}
\newcommand{\cate}[1]{\mathcal{#1}}
\newcommand{\framme}[1]{\mathfrak{#1}}
\newcommand{\divi}[1]{\mathrm{div}(#1)}
\renewcommand{\diamond}{\lozenge}
\newcommand{\card}[1]{\vert #1 \vert}
\newcommand{\regle}[1]{(\mathrm{#1})}

\newtheorem{thm}{Theorem}[section]

\newtheorem{lem}[thm]{Lemma}
\newtheorem{prop}[thm]{Proposition}
\theoremstyle{definition}
\newtheorem{defn}[thm]{Definition}
\theoremstyle{remark}

\numberwithin{equation}{section}

\title[Compl. results for MV  \L ukasiewicz modal syst. and rel. semant.]{Completeness results for many-valued \L ukasiewicz modal systems and relational semantics}
\author[G. Hansoul and B. Teheux]{Georges Hansoul and Bruno Teheux}
\email{g.hansoul@ulg.ac.be, b.teheux@ulg.ac.be}
\address{D\' epartement de Math\' ematiques, Universit\' e de Li\` ege, 12,
  Grande Traverse, 4000 Li\` ege, Belgium.}
\begin{document}
\begin{abstract}
The paper is dedicated to the problem of adding a modality to the \L
ukasiewicz many-valued logics in the purpose of obtaining completeness results
for \textsc{Kripke} semantics. We define a class of modal many-valued logics
and their corresponding \textsc{Kripke} models and modal many-valued algebras. Completeness results are
considered through the construction of a canonical model. Completeness is
obtained for modal finitely-valued logics but also for a modal many-valued
system with an infinitary deduction rule. We introduce two classes of frames for
the finitely-valued logics and show that they define two distinct classes of
\textsc{Kripke}-complete logics.
\end{abstract}
\maketitle


\section{Introduction}
Modal logics and many-valued logics were both historically introduced in order
to free oneself from the rigidity of propositional logic. With many-valued
logics, the logician can choose the truth values of the propositions in a set
with more than two elements. With modal logics, the logician introduce a new
connector whose aim is, for instance, to model the possibility. Even if these two
approaches of the science of reasoning are not born  independently (see
chapter 21 of \cite{Got}),
many-valued logics and modal logics received distinguished interests since
their birth.

On the one hand, mathematicians tackled many-valued logics 
(as defined by \textsc{J. \L ukasiewicz} in \cite{Luka1}; see \cite{Luka2} for
an English translation and \cite{Cign} for a monograph on the subject)
through their algebraic form: the class of MV-algebras that was introduced by
\textsc{C.C. Chang} in 1958 (see \cite{Chang1} and \cite{Chang2}) in order to
obtain an algebraic proof of the completeness result for the
infinite-valued \textsc{\L ukasiewicz} logic. 

On the other hand, modal logics were also  studied  through their algebraic
disguises, which are the Boolean algebras with operators (introduced in
\cite{Tar1} and \cite{Tar2}). But the success of modal logics among the communities of mathematicians, computer
scientists and philosophers is a consequence of the relational semantics
introduced in the sixties by \textsc{S. Kripke} (see \cite{Kripke}). With  \textsc{Kripke}
semantic, also called possible worlds semantic, a formula is possible in a
world $w$ if it is true in a world accessible from $w$.  From then on, in
their approach of modal logic, mathematicians have been focusing their
attention on the connection between the algebraic and the relational
semantics. This approach allowed a great simplification of the proof of the
completeness of normal modal logics through the construction of the so-called
canonical model (see \cite{Vene1} for example).  

Since the definition of a \textsc{Kripke} model can easily be extended to a
many-valued realm, the problem of merging modal and many-valued logic has already been considered
by some mathematicians in the last few years (see \cite{Fit1} and \cite{Fit2} or \cite{Ost} for example). But the algebraic
approach and its connection with the relational semantics has never been
considered.

In this paper, we introduce some modal many-valued (in a \textsc{\L
ukasiewicz} meaning) logics and their corresponding algebras and tackle the
problem of the completeness with respect to many-valued \textsc{Kripke}
models through the construction of a canonical model.

The next section of this paper is dedicated to the introduction of the many-valued
\textsc{Kripke} models. In the third section, we define the modal many-valued
logics. The axiomatization of these logics is very natural since it appears
clearly that they admit the class of many-valued \textsc{Kripke} models as a
sound semantic. 

Modal many-valued algebras are introduced in the fourth section
as a step toward the construction of the canonical model of a modal
many-valued logic, which is the object of the next section.  Here, the reader
can find one of the main result of this paper: the natural definition of the
valuation on the canonical model extends to
formulas (Proposition \ref{prop:CanoMod}). 

Completeness results with respect to the \textsc{Kripke} models of section
\ref{sect:RelSem} are considered in section  \ref{sect:compl}. The results are
obtained for the finitely-valued logics but also for an infinitary modal
many-valued system (i.e. a formal system that admits an infinitary deduction
rule). The construction of the canonical model allows us to
simplify the axiomatization of the modal finitely-valued logics. We also
introduce the class of $n+1$-frames. Roughly speaking, an $n+1$-frame is
a first-order structure obtained from a frame by deciding to ban some valuations on the frame.  We then
illustrate the dissimilarity between \textsc{Kripke}-completeness and
$n+1$-\textsc{Kripke}-completeness. 

We conclude the paper by suggesting some tracks for possible projects.

\section{A relational semantic for $[0,1]$-valued modal logics}\label{sect:RelSem}

The modal many-valued systems that we develop in the sequel of the paper admit
a sound (and for some of them complete) relational semantic. Since this
semantic is the main strength of these systems, we have decided to first
introduce it.


Let us denote by $\Prop$ an infinite set of  propositional variables and by $\Form$ the set of formulas defined inductively by the following rules:
\begin{enumerate}
\item $\Prop \subseteq \Form$;
\item if $\phi$ and $\psi$ are in $\Form$ then $\neg \phi$, $\phi \oplus \psi$ and $\square \phi$ are in $\Form$.
\end{enumerate}
The intended meaning of $\phi \oplus \psi$ and $\neg \psi$ is clear (these formulas have their usual \textsc{\L ukasiewicz} meaning) and $\square \phi$ can be read, for example, as ``$\phi$ necessary holds''.

Obviously, as usual, we make use of the following abbreviations: the formula
$\phi \implies \psi$ stands for $\psi \oplus \neg \phi$, the formula $\psi
\odot \phi$  for $\neg (\neg \psi \oplus \neg \phi)$, the formula $\phi \vee \psi$ for $(\phi
\implies \psi) \implies \psi$, the formula $\phi \wedge \psi$ for $\neg(\neg \phi
\vee \neg \psi)$ and $\diamond \phi$ for $\neg \square \neg \phi$.

The definition of a \textsc{Kripke} model can easily be extended to a
$[0,1]$-valued realm. In the following, we consider the real
unit interval $[0,1]$ endowed with the \textsc{\L ukasiewicz} implication and negation: 
$x \implies y= \min (1, 1-x+y) \mbox{ and } \neg x=1-x$. If $n$ is a positive
integer, we denote by $\lucas_n$ the subset $\{0, \frac{1}{n},
\ldots, \frac{n-1}{n},1\}$ of $[0,1]$. Note that $\lucas_n$ is closed for $\implies$ and $\neg$.
\begin{defn}
A \emph{many-valued \textsc{Kripke} model} (or simply a \emph{many-valued model}) $\model{M}=\struc{W, R, \Val}$ is
given by a non empty set $W$, an accessibility relation $R \subseteq W \times
W$ and a map $\Val: \Prop \times W \rightarrow [0,1]$. If $n$ is a
positive integer such that $\Val(\Prop,
W)\subseteq \lucas_n$, then $\model{M}$ is called an \emph{$n+1$-valued
  \textsc{Kripke} model}.

A \emph{frame} $\framme{F}=\struc{W,R}$ is given by an non empty set $W$ and an
accessibility relation $R$ on $W$. A model $\model{M}=\struc{W',R', \Val}$ is
\emph{based on the frame $\framme{F}=\struc{W,R}$} if $W=W'$ and $R=R'$.
\end{defn}
If $\model{M}=\struc{W, R, \Val}$ is a many-valued \textsc{Kripke} model, we
extend inductively the map $\Val$ to formulas of $\Form$ by the following
rules:
\begin{itemize}
\item $\Val(\phi \oplus \psi, w)=\Val(\phi, w) \oplus \Val(\psi, w)$ and
  $\Val(\neg \phi, w)=\neg \Val(\phi, w)$,
\item $\Val(\square \phi, w)= \bigwedge \{\Val(\phi, w') \mid (w,w') \in R\}$,
\end{itemize}
for any formulas $\phi$ and $\psi$ of $\Form$ and any world $w$ of $W$ (where
$\bigwedge$ means that we consider the infemum in $[0,1]$). 

We write  $\model{M}, w \models \phi$ or simply $w \models \phi$ (and
say that $w$ \emph{satisfies} $\phi$)
whenever $\Val(\phi, w)=1$, and $\model{M} \models
\phi$ whenever $w \models \phi$ for any $w$ in $W$. In that case, we say
$\phi$ is \emph{true} in $\model{M}$. Formulas that are true in any model
$\model{M}$ are called \emph{tautologies}. If $\Gamma$ is a set of formulas
that are true in a model $\model{M}$, then $\model{M}$ is a \emph{model of
  $\Gamma$}. If $\framme{F}$ is a frame and $\phi$ is a formula that is true
in any model based on $\framme{F}$, we say that $\phi$ is \emph{valid} in
$\framme{F}$ and write $\framme{F}\models \phi$.

Note that the $2$-valued \textsc{Kripke} models coincide with the
\textsc{Kripke} models of normal modal logics (where the operation
$\oplus$ matches up with the supremum $\vee$).

\begin{prop}\label{prop:tauto}
  If $\tau$ is an increasing unary term of the language
  $\lang_{MV}=\{\implies, \neg\}$, then the formulas
\[
\begin{array}{ccc}
\square (p \implies q)\implies (\square p\implies \square q), &  \square (p \wedge q)\leftrightarrow\square p \wedge
\square q, & \square \tau(p) \leftrightarrow \tau(\square p)
\end{array}\]
are tautologies.
\end{prop}
The idea of using these models as a semantic for modal many-valued  logics is
not new. See \cite{Ost} for example.

\section{\L ukasiewicz modal many-valued logics}
The purpose of this section is to introduce a family of \emph{modal many-valued 
  logics} and their corresponding algebras in order to tackle completeness
  results through the construction of a \emph{canonical model}. We refer to
  \cite{Cign} for an introduction to \textsc{\L ukasiewicz} logic and to
  \cite{Vene1} and \cite{Chagrov} for an introduction to modal logic. 

  \begin{defn}\label{def:mmvl}
    A \emph{modal many-valued logic} is a set $\logic{L}$ of formulas of
    $\Form$ that is  closed under \emph{modus ponens}, substitution, the
    necessitation rule $\regle{RN}$ (from $\phi$ infer $\square \phi$)   and
    that contains
    \begin{itemize}
    \item an axiomatic base of \textsc{\L ukasiewicz} many-valued logic ($p \implies (q \implies
   p)$, $(p \implies q)\implies ((q \implies r)\implies(p
   \implies r))$, $((p \implies q)\implies q) \implies ((q
   \implies p)\implies p)$, $(\neg p \implies \neg q)\implies
   (q \implies p)$ for example);
    \item the axiom $(\K)$ of modal logic: $\square(p \implies q)
   \implies (\square p \implies \square q)$,
    \item the formulas $\square (p \oplus p) \leftrightarrow \square p \oplus
    \square p$ and $\square (p \odot p) \leftrightarrow \square p \odot
    \square p$,
    \item the formula $\square(p \oplus p^m) \leftrightarrow (\square p \oplus (\square
    p)^m)$ for every positive integer $m$.
    \end{itemize}
 As usual, we  write $\proves_\logic{L} \phi$ and say that $\phi$ is a
\emph{theorem of $\logic{L}$} whenever $\phi \in \logic{L}$ and denote by
$\logic{K}$ the smallest modal many-valued logic. If in
addition $\logic{L}$ contains an axiomatic base of the $n+1$-valued \textsc{\L
ukasiewicz} logic, we say that $\logic{L}$ is a \emph{modal $n+1$-valued  
logic} and we denote by $\logic{K}_n$ the smallest of these logics. 
\end{defn}
 Note that, according to Proposition
\ref{prop:tauto}, the proposed axioms are tautologies for the many-valued \textsc{Kripke}
models, so that the latters form a sound semantic for the modal many-valued
logics.

 Let us also remark that, as it will appear clearly in the sequel of the paper (in
 Proposition \ref{prop:CanoMod}), we only use
 the last family of axioms as a kind of conservative law for $\square$ with respect
 to infinitely great elements. Moreover Proposition \ref{prop:simpliaxio}
 gives an axiomatization of the finitely-valued logics without this family of
 axioms (and this explains why we have added the axiom $\square (p \oplus p) \leftrightarrow \square p \oplus
    \square p$  even if it is equivalent to $\square(p \oplus p^m) \leftrightarrow (\square p \oplus (\square
    p)^m)$ with $m=1$).  

We can easily gather the following theorems and admissible rules of
$\logic{K}$.

\begin{prop}
The following formulas are theorems of $\logic{K}$:
\begin{center}\begin{tabular}{ll}
 $\square(p \implies q) \implies (\diamond p \implies \diamond q)$, & $\diamond (p \oplus q) \implies (\diamond p \oplus \diamond q)$,\\
 $(\square p \wedge \diamond q) \implies \diamond (p \wedge q)$, & $\square (p
 \wedge q) \implies (\square p \wedge \square q)$ \\ $(\square p \odot \square q) \implies
 \square (p \odot  q)$.
\end{tabular}
\end{center}
Moreover, the logic $\logic{K}$ is  closed under the following deduction
rules:
\begin{enumerate}
 \item $
    \begin{array}{c}
     \phi \leftrightarrow \psi \\ \hline
     \square \phi \leftrightarrow \square \psi
    \end{array}$,
  \item  $\begin{array}{c}
     \phi_1 \odot \cdots \odot \phi_n  \rightarrow \psi \\ \hline
    \square \phi_1 \odot \cdots \odot \square \phi_n  \rightarrow \square \psi  
    \end{array}$.
\end{enumerate}
\end{prop}
\begin{proof}
  The proofs are simple adaptations of the two-valued proofs. 
\end{proof}
Note that at this point of our development, we can not decide if the 
formula $\square p \wedge  \square q \rightarrow \square (p \wedge  q)$ is in
$\logic{K}$ or in $\logic{K}_n$ (for $n \geq 2$). We shall conclude latter, thanks to a completeness
result, that it is a theorem of $\logic{K}_n$ for any $n$. 
On the opposite, the formula  $\square (p \odot q) \implies \square p \odot
\square q$ is not a theorem of $\logic{K}_n$ for any $n\geq 2$ since it is not a $\logic{K}_n$-tautology.
\begin{defn}
If $\Gamma \cup \{\phi\}$ is a set of formulas and if $\logic{L}$ is a many-valued
modal logic, we say that \emph{$\phi$ is deducible from $\Gamma$ in
  $\logic{L}$} and write $\Gamma \proves_{\logic{L}} \phi$ (or simply
$\Gamma \proves \phi$ when $\logic{L}=\logic{K}$ or $\logic{L}=\logic{K}_n$
following the context) if $\phi$ is in any extension of $\Gamma \cup \logic{L}$
that is closed under substitution and \emph{modus ponens}.
\end{defn}
Note that we can state the following adaptation of the deduction theorem.
\begin{lem}
  If $\Gamma \cup \{\phi\}$ is a set of formulas and if $\logic{L}$ is a
  modal many-valued logic, then $\Gamma \proves_{\logic{L}} \phi$ if and only
  if there is a finite subset $\{\phi_1, \ldots, \phi_r\}$ of $\Gamma$ and some
  positive integers $m_1, \ldots, m_r$  such that $\proves_\logic{L}
  \phi_1^{m_1} \odot \cdots \odot \phi_r^{m_r}$.
\end{lem}


\section{Modal many-valued  algebras and the algebraic semantic}\label{sect:logic}
We introduce very briefly the varieties of \emph{modal many-valued  algebras}
and state the completeness result for modal many-valued logics  and
algebras. This somehow obvious result can be seen as a step towards the
construction of the \emph{canonical model} and a possible completeness theorem for
many-valued \textsc{Kripke} models. We refer the reader to \cite{Cign} or
\cite{MvAlgUniv} for an introduction to the variety of MV-algebras. 

\begin{defn}
  If $\logic{L}$ is a modal many-valued  logic then an
  \emph{$\logic{L}$-algebra} is an algebra $A$ over the language
  $\lang_{MMV}=\{\implies, \neg, \square, 0, 1\}$ that satisfies the equations
  naturally induced by the formulas of $\logic{L}$. We denote by $\var{MMV}$
  (resp. $\var{MMV}_n$) the variety of $\logic{K}$-algebras (resp. the variety
  of $\logic{K}_n$-algebras). Members of $\var{MMV}$ (resp. $\var{MMV}_n$) are simply called
  \emph{modal many-valued  algebras} or MMV-algebras
  (resp. \emph{modal $n+1$-valued  algebras} or MMV$\phantom{}_n$-algebras).
\end{defn}
A modal many-valued logic $\logic{L}$ is often given by a set $\Gamma$ of axioms,
i.e. $\logic{L}$ is the smallest modal many-valued logic that contains
$\logic{K} \cup \Gamma$, and is denoted by $\logic{K} + \Gamma$.

Since the most commonly used axiomatization of the variety $\var{MV}$ of MV-algebras is given over the
language $\{\oplus, \odot, \neg, 0, 1\}$, we  preferably  use this language
instead of $\lang_{MV}$ (with the help of the theorem $(\phi \oplus \psi)
\leftrightarrow (\neg \phi \implies \psi)$). Thus, an MMV-algebra is
simply an algebra $A=\struc{A, \oplus, \odot, \neg, \square, 0, 1}$ of type
$(2,2,1,1,0,0)$ such that 
\begin{itemize}
\item the reduct of $A$ to the language $\{\oplus, \odot, \neg, 0, 1\}$ is an
  MV-algebra (i.e. A satisfies the equations $\neg \neg
x=x$, $x \oplus 1=1$, $\neg 0=1$, $x \odot y=\neg(\neg x \oplus \neg y)$,
$(x \odot \neg y) \oplus y=(y \odot \neg x) \oplus x$);
\item the algebra $A$ satisfies the equations $\square (x \implies y) \implies
  \square  x \implies \square y=1$, $\square (x \oplus x) = \square x \oplus \square x$,
  $\square (x \odot x)= \square x \odot \square x$ and $\square (x
  \oplus x^m) = (\square x \oplus (\square x)^m)$ for every positive
  integer $m$.
\end{itemize}
Similarly, an MMV${}_n$-algebra is an MMV-algebra whose
reduct to the language of MV-algebras is a member of the variety
$\classop{HSP}(\lucas_n)=\classop{ISP}(\lucas_n).$

Recall that on an {MV}-algebra $A$, the relation $\leq$ defined by 
\[
x \leq y \mbox{ if } x \implies y=1
\]
is a lattice order on $A$ with $x \vee y=(x \implies y) \implies y$
and $x \wedge y=\neg (\neg x \vee \neg y)$.

It is not the object of this paper to study the varieties of
{MMV}-algebras in details. This work should be done in a forthcoming
paper. Let us note that the duality developed in \cite{BruDual} for each of
the varieties $\mathcal{MMV}_n$ could be used as a tool for the investigation of
these varieties.
\begin{defn}
A \emph{filter} of an MMV-algebra $A$ is a filter of its
MV-algebra reduct (i.e. a non empty increasing subset of $A$ that
contains $y$ whenever it contains $x$ and $x \implies y$). If $x$ is a subset
of $A$, we denote by $\struc{X}$ the filter generated by $X$.

We denote by $\free_\logic{L}$ the free $\logic{L}$-algebra over the set
$\Prop$ of variables for any modal many-valued logic $\logic{L}$, i.e. the
set of formulas modulo $\logic{L}$-equivalence. In the following, we do
not distinguish a formula $\phi$ with its  class $\phi^{\logic{L}}$  in
$\free_\logic{L}$. 
\end{defn}

Recall that the lattice of filters of an MV-algebra $A$ is isomorphic to
the lattice of congruences of $A$. The congruence $\theta_F$ associated to a filter $F$
by this isomorphism is defined by $(x,y) \in \theta_F$ if $(x \implies y)
\odot (y \implies x) \in F$. As usual, we denote by $A/F$ the quotient $A/\theta_F$.

 For our purpose, the next result is fundamental, albeit an obvious
one. 
\begin{prop}
  If $\logic{L}$ is a modal many-valued logic, and $\Gamma \cup \{\phi\}$ is a
  set of formulas then $\Gamma \proves_\logic{L} \phi$ if and only if $\phi=1$
   in $\free_{\logic{L}}/\struc{\Gamma}$ or equivalently if $\phi$ is satisfied
  in every $\logic{L}$-algebra  that satisfies the axioms of
  $\Gamma$.
\end{prop}


\section{Construction of the canonical model}
Here is one of the main contributions of the paper. Recall that the variety of
 {MV}-algebras
 is the variety generated by the algebra $\struc{[0,1], \oplus,
  \odot, \neg, 0, 1}$ defined on the real unit interval $[0,1]$ by $x \oplus
y=\min(x+y,1)$ and $\neg x= 1-x$ and that an {MV}-algebra $A$ is
simple if and only if it is an isomorphic copy of a subalgebra of
 $[0,1]$. Moreover,  two isomorphic subalgebras of $[0,1]$ are necessarily
 equal (and the isomorphism is the identity). We can thus state the following
 lemma, which will enable us to define a valuation on the canonical model. A
 filter of an MV-algebra $A$ is \emph{maximal} if it is maximal among the
 proper filters of $A$.
 
\begin{lem}\label{lem:filtermax}
A filter $F$ of an MV-algebra $A$ is maximal  if and only if there is a unique
embedding $f:A/F \rightarrow [0,1]$.
\end{lem}

The idea of the construction of the canonical model for a logic $\logic{L}$ is
classical. The universe of the model is the set of the maximal filters of
$\free_L$ (it coincides with the maximal consistent extensions of
$\logic{L}$). But, in order to simplify the expression of our results, it is
better to identify, with the help of Lemma \ref{lem:filtermax}, the set of the maximal filters of $\free_\logic{L}$ with
the set $\var{MV}(\free_\logic{L}, [0,1])$ of the homomorphisms of MV-algebras
from $\free_\logic{L}$ to $[0,1]$.

\begin{defn}
The \emph{canonical model} for a modal many-valued logic $\logic{L}$ is the
model $\model{M}_\logic{L}=(W_\logic{L}, R_\logic{L}, \Val_{\logic{L}})$ where
\begin{itemize}
\item the universe $W_\logic{L}$ of $\model{M}_\logic{L}$ is the set
  $\var{MV}(\free_\logic{L}, [0,1])$;
\item the binary relation $R_\logic{L}$ is defined by
\[
(u,v) \in R_\logic{L} \mbox{ if } \forall \ \phi \in \free_L \ (u(\square \phi)=1
\Rightarrow v(\phi)=1), 
\]
\item the valuation map $\Val_\logic{L}: \Prop \times W_\logic{L} \rightarrow
  [0,1]$ is defined by \[Val_\logic{L}(u,p)=u(p).\]  
\end{itemize}  
\end{defn}
Note that the definition of the canonical model for $\logic{K}_1$ coincides with the classical
definition of the canonical model for the Boolean basic modal logic (if we
identify any maximal filter $F$ of $\free_{\logic{K_1}}$ with the quotient map
$\pi_F: \free_{\logic{K_1}} \rightarrow \{0, 1\}$ and if we identify the
Boolean valuation map $\Val: \Prop \rightarrow
\mathcal{P}(W_{\logic{K}_1})$ with its characteristic  function).

The main result of this section is that the map $\Val_{\logic{L}}$ extends to
formulas. Before considering the proof of this result, we need the following
definition.
\begin{defn}\label{defi:termes}
We denote by $\2$ the subset of $\Q$ that contains the numbers that can be
written as a finite sum of powers of 2. If $r$ is in $\2 \cap [0,1]$, we denote by $\tau_r$ a composition of the terms
$x \oplus x$ and $x \odot x$ such that $\tau_r(x) < 1$  for every  $x \in
[0,r[$  and $\tau_r (x) = 1$ for every $x \in [r,1]$. A proof of the
existence of such terms can be found in \cite{Ost} for example.
Furthermore, we can always choose $\tau_r$ such that $\tau_r (x)=1$ for
every $x \in \lucas_n \cap [r,1]$ (but this choice is not independant of
$n$).
\end{defn}
\begin{lem}\label{lem:Requiv}
  If $\logic{L}$ is a modal many-valued logic and if $u, v \in W_\logic{L}$, then $(u,v) \in R_\logic{L}$ if and only if $u \circ
  \square \leq v$.
\end{lem}
\begin{proof}
The right to left part of the assertion is clear. Let us prove
  the left to right part and suppose that there is a $\phi$ in $\free_\logic{L}$, a
  $v$ in $R_\logic{L} u$ and
  an $r$ in $\2 \cap [0,1]$ such
  that $v(\phi) < r \leq u(\square \phi)$. It follows that
\[
\begin{array}{ccc}
\tau_r(v(\phi))=v(\tau_r(\phi)) <1 & \mbox { and } & 1=\tau_r(u(\square
\phi))=u(\tau_r(\square \phi))=u(\square \tau_r(\phi)),
\end{array}
\]
which is a contradiction since $u R_{\logic{L}}v$.
\end{proof}
\begin{prop}\label{prop:CanoMod}
  If $\logic{L}$ is a modal many-valued logic, then
\[
\Val_{\model{M}_{\logic{L}}}(\phi,u)=u(\phi)
\]
for any $\phi$ in $\Form$ and $u$ in $W_\logic{L}$.
\end{prop}
\begin{proof}
The non trivial statement is the equality
\begin{equation}\label{eqn:these}
u(\square \phi)=\bigwedge_{v \in R_{\logic{L}}u}  v(\phi).
\end{equation}
The inequality $\leq$ is the content of Lemma \ref{lem:Requiv}.

Suppose now that  the equality does not hold in (\ref{eqn:these}), but just
the strict inequality $<$. Then, there is an $r$ in $\2 \cap [0,1]$ such that
\[
u(\square \phi) < r \leq \bigwedge_{v \in R_{\logic{L}}u}  v(\phi).
\]
It means that for any $v \in R_\logic{L}u$, the maximal filter $v^{-1}(1)$ of
$\free_\logic{L}$  contains $\tau_r(\phi)$ and that the filter $\square^{-1}
u^{-1}(1)$ does not contain $\tau_r(\phi)$. 

But, since $(u,v) \in R_\logic{L}$ if and only if  $\square^{-1}u^{-1}(1)
\subseteq v^{-1}(1)$, if follows that the maximal filters of $\free_\logic{L}$
that contain $\square^{-1}u^{-1}(1)$ are exactly the $v^{-1}(1)$ with $v \in
R_\logic{L} u$, while each of these maximal filters contains $\tau_r(\phi)$. It
means that the class of $\tau_r(\phi)$ in $\free_\logic{L} / \square^{-1}
u^{-1}(1)$ is infinitely great, so that $\tau_r(\phi) \oplus \tau_r(\phi)^m$
belongs to $\square^{-1} u^{-1}(1)$ for every positive integer $m$. It follows
that
\[1=u(\square (\tau_r(\phi) \oplus \tau_r(\phi)^m)) \leq u(\tau_r(\square
\phi) \oplus (\tau_r(\square \phi))^m), \]
for any positive integer $m$, so that $u(\tau_r(\square \phi))$ is infinitely
great in $u(\free_\logic{L})$. Since $u(\free_\logic{L})$ is a subalgebra of $[0,1]$, we
obtain that $u(\tau_r(\square \phi))=1$, a contradiction.
\end{proof}


\section{Completeness results}\label{sect:compl}

Proposition \ref{prop:CanoMod} enables to prove completeness results
for some modal many-valued logics\ldots~but only for some of them. Indeed,
with the help of Proposition \ref{prop:CanoMod}, we obtain that if $\Gamma$ is
a set of axioms then a formula $\phi$ that is valid in every model of $\Gamma$
is deducible from $\Gamma$ if the algebra $\free_\logic{K}/\struc{\Gamma}$ is
semi-simple, i.e. if $\struc{\Gamma}$ is the only element of $\free_\logic{L}/\struc{\Gamma}$ that
is in every maximal filter of $\free_\logic{L}/\struc{\Gamma}$. Unfortunately,
contrary to the two-valued case (where every logic is the intersection of its
maximal consistent extensions), we can not ensure  \emph{a priori} that
$\free_\logic{L}/\struc{\Gamma}$ is a semi-simple algebra.

Hopefully, there are some very interesting logics for which the completeness
result can be stated.
\subsection{Modal finitely-valued logics}
The first family of systems that admit the many-valued \textsc{Kripke} models
as a complete semantic is the finitely-valued ones. 
\begin{thm}
  If $\Gamma\cup \{\phi\}$ is a set of formulas, then $\Gamma
  \proves_{\logic{K}_n} \phi$ if and only if $\phi$ is valid in every
  $n+1$-valued \textsc{Kripke} model of $\Gamma$.
\end{thm}
\begin{proof}
  The algebra $\free_{\logic{K}_n}/\struc{\Gamma}$ is a member of
  $\classop{HSP}(\lucas_n)=\classop{ISP}(\lucas_n)$ and so is semi-simple. 
\end{proof}
Note that by considering $n=1$, the  preceding proposition boils down to the
completeness result for  Boolean  basic modal logic and \textsc{Kripke}
semantic.

We have announced in section \ref{sect:logic} the following result which is an
application of the preceding completeness theorem.
\begin{prop}
 If $n$ is a positive integer then $\proves_{\logic{K}_n} \square (p \wedge q ) \leftrightarrow (\square p \wedge \square
  q)$. 
\end{prop}

Moreover, we can simplify the axiomatization \ref{def:mmvl} of
$\logic{K}_n$. We can indeed get rid off the family of axioms that expresses the
conservative law of $\square$ with respect to the infinitely great elements. 
\begin{prop}\label{prop:simpliaxio}
  If $\logic{MV}_n$ denotes  the $n+1$-valued \textsc{\L ukasiewicz}
  logic and if $\logic{K'}_n=\logic{MV}_n + \square (p \implies q) \implies
  (\square p \implies \square q) + \square (p \oplus p) \leftrightarrow
  (\square p \oplus \square p) +\square (p \odot p) \leftrightarrow   (\square
  p \odot \square p)$, then $\logic{K}_n=\logic{K'}_n$.
\end{prop}
\begin{proof}
In the proof of Proposition \ref{prop:CanoMod} with $\logic{L}=\logic{K}_n$,
we can deduce directly  that $u(\square  \tau_r
(\phi))=1$ from the fact that $\tau_r(\phi)$ is infinitely great
in $\free_{\logic{K}_n}/\square^{-1}u^{-1}(1)$, since $\free_{\logic{K}_n}/\square^{-1}u^{-1}(1)$ has no non trivial
infinitely great element. It means that Proposition \ref{prop:CanoMod} stands
with $\logic{L}=\logic{K}_n'$ and that $\proves_{\logic{K}_n'} \phi$ for any
formula $\phi$ that is valid in any $n+1$-valued \textsc{Kripke} model. We can
thus conclude since for any positive integer $m$, the formula $\square (p
\oplus p^m) \implies (\square p \oplus (\square p)^m)$ is a tautology.
\end{proof}

Apart from the completeness result, the extensions of $\logic{K}_n$ seem to
share interesting properties with the Boolean modal logics. For instance, the
paper \cite{BruDual} is dedicated to the construction of a duality for
$\logic{K}_n$-algebras and a class of topological structures. For $n=1$, this
duality coincides with the \textsc{Stone} duality for modal algebras.
The role played by the duality for $\logic{K}_n$-algebras is as important  as
the role played by the \textsc{Stone} duality for modal algebras in Boolean modal
logic. Indeed, the class of the dual structures (called $\cate{MX}_n$-structures) forms a very adequate semantic
since any extension $\logic{L}$ of $\logic{K}_n$ is complete with respect to
the class of $\cate{MX}_n$-structures in which $\logic{L}$ is valid.

Moreover, the construction of this duality suggests two ways of going from
$n+1$-valued \textsc{Kripke} models to frames (and conversely). Indeed, we can obviously
define a frame $\framme{F}=\struc{W, R}$ as a a set $W$ with a binary relation $R$
on $W$. Then, a frame $\framme{F}=\struc{W, R}$ becomes a model by the
addition of a valuation $\Val:\Prop \times W \rightarrow \lucas_n$. So,
the set of truth values in a world $w$ is given by the valuation, at the model
level and not at the frame level.  
Now, it is also possible (and as we shall see, relevant) to consider some new
(first order) structures, called \emph{$n+1$-frames}, in which the set of truth values in a world $w$ is know \emph{a
  priori}, without any reference to a valuation. In the following definition,
we denote by $\divi{n}$ the set of the positive divisors of $n$.
\begin{defn}
  An \emph{$n+1$-frame} $\framme{F}=\struc{W, \{r_m \mid m \in \divi{n}\}, R}$
  is given by a set $W$, a subset $r_m$ of $W$ for every $m$ in $\divi{n}$ and a
  relation $R \subseteq W \times W$ such that
  \begin{enumerate}
    \item for every $m$ and $k$ in $\divi{n}$, the intersection $r_m \cap r_k$
    coincides with $r_{\gcd (m,k)}$ and $r_n=W$;
    \item for every $m$ in $\divi{n}$, the set $Rr_m=\{w \mid \exists \ w' \in
    r_m \ w'Rw\}$ of the successors of the elements of $r_m$ is a subset of $r_m$.
  \end{enumerate}

A model $\struc{W', R', \Val}$ is \emph{based on an the $n+1$-frame}
$\framme{F}=\struc{W, \{r_m \mid m \in \divi{n}\}, R}$ if $W=W'$, $R=R'$ and
$\Val(p, w) \in \lucas_m$ for any $m$ in $\divi{n}$, any $w$ in $r_m$ and $p$
in $\Prop$.
\end{defn}

Validity in $n+1$-frames is defined similarly as in the class of frames. 

Thus, an $n+1$-frame is obtained from a frame by
restricting the class of valuations that can be added to this frame to define
an $n+1$-\textsc{Kripke} model.
We should so have a gain in the expressivity of the class of
$n+1$-frames with regard to the class of frames. For instance, there are
some extensions $\logic{L}$ of $\logic{K}_n$ that are characterized by a class
of $n+1$-frames but that are not characterized by any class of frames. Here
are a few easy examples.
\begin{defn}
  A modal many-valued logic $\logic{L}$ is \emph{\textsc{Kripke} complete}
  (resp. \emph{tabular}) if there is
  a class of frames $K$ (resp. a finite frame $\framme{F}$) such that $\logic{L}$ is the set of formulas
  that are valid in every frame of $K$ (resp. in $\framme{F}$). 

  Similarly, a logic $\logic{L}$ is  \emph{$n+1$-\textsc{Kripke} complete}
  (resp. \emph{$n+1$-tabular}) if there
  is a class of $n+1$-frames (resp. an $n+1$-frame $\framme{F}$) such that $\logic{L}$ is the set of formulas that
  are valid in every model based on a frame of $K$ (resp. on $\framme{F}$).

  If $\framme{F}=\struc{W,\{r_m \mid m \in \divi{n}\}, R}$ and
  $\framme{F}=\struc{W',\{r'_m \mid m \in \divi{n}\}, R'}$ are two
  $n+1$-frames, a map $f:W \rightarrow W'$ is called a \emph{$n+1$-$\pi$-morphism}
  if the three following conditions are satisfied:
  \begin{enumerate}
  \item if $u$ and $v$ are in $W$ and $(u,v) \in R$ then $(f(u),f(v)) \in R'$;
  \item if $u \in W$ and $v'\in W'$ with $(f(u), v') \in
  R$ then there exists a $v$ in $Ru$ such that $f(v)=v'$;
  \item if $u \in r_m$ then $f(u) \in r'_m$.
  \end{enumerate}
\end{defn}
We leave to the reader the task to prove that validity is preserved under
$n+1$-$\pi$-morphic image, i.e. that if $f:\framme{F} \rightarrow \framme{F}'$
is a surjective $n+1$-$\pi$-morphism between two $n+1$-frames $\framme{F}$ and
$\framme{F}'$, then $\framme{F'} \models \phi$ whenever $\framme{F}\models \phi$.
\begin{prop} Assume that $n \geq 2$. We have the following completeness results.
  \begin{enumerate}
  \item The logic $\logic{L}_1= \logic{K}_n+(\square p \vee \square \neg p)$ is $n+1$-\textsc{Kripke} complete with respect
  to the class of the $n+1$-frames that satisfy $\forall u \ Ru \subseteq
  r_1$ but is not \textsc{Kripke}-complete.
  \item The logic $\logic{L}_2=\logic{K}_n+\square (p \vee \neg p)+\square
  \square p + \neg (\diamond p \wedge \diamond \neg p)$ is
  $n+1$-tabular but is not even \textsc{Kripke} complete.
  \item The logic $\logic{K}_n + \square x \implies x$ is \textsc{Kripke}-complete with respect to
  the class of reflexive frames.
  \item The logic $\logic{K}_n + \square x \implies \square \square x$ is
  \textsc{Kripke}-complete with respect to the class of transitive frames.
  \end{enumerate}
\end{prop}
\begin{proof}
  (1) It is clear, by the definition of a model based on an $n+1$-frame that
      $\logic{L}_1$ is characterized by the class of $n+1$-frames that
      satisfy the first order formula  $\forall u \ Ru \subseteq
      r_1$.

      Suppose then that $K$ is a class of frames that characterizes
      $\logic{L}_1$. First note that we may suppose that $K$ contains a frame
      $\struc{W,R}$ with $R$ non trivial. Otherwise, for any frame $\framme{F}$
      of $K$ and any formula $\phi$, the formula $\square \phi$ is valid in
      $\framme{F}$ and is so a theorem of $\logic{L}_1$. But it is easy to
      construct a $\logic{L}_1$-counter-model for $\square (p \wedge \neg p)$.

      Now, if $\framme{F}=\struc{W,R}$ is a frame of $K$ with a  non trivial
      relation $R$ and
      if $\model{M}=\struc{W,R, \Val}$ is a model based on $\framme{F}$ and
      $w, v \in W$ with $wRv$, it follows that $\model{M}, w \models (\square p \vee \square
      \neg p)$. We deduce  that $\Val(p, v) \in \{0,1\}$. Then, if we denote by $\model{M}'=\struc{W,R,\Val'}$
      the model based on $\framme{F}$ defined by 
      \[\Val'(q, u)=\left\{\begin{array}{l}\Val(q,u) \mbox{ if } q \neq p \mbox{
      or } u \neq v,\\ \frac{1}{n} \mbox{ if } q=p \mbox{ and }
      u=v,\end{array}\right.\]
      it appears that $\square p \vee \square \neg p$ is not true in
      $\model{M}'$, a contradiction since $\model{M}'$ is based on a frame of $K$.
  
  (2) Let us consider the $n+1$-frame $\framme{F}$ whose universe is
      $\{u,v\}$ with $(u,v) \in R$, $u \in r_n$ and $v \in r_1$ (we only
      specify for any world of an $n+1$-frame the smallest of the subsets
      $r_m$ that contain this world) and the $n+1$-frame $\framme{F}'$ as the one
      irreflexive point belonging to $r_n$.

      It is clear that any formula of $\logic{L}_2$ is satisfied in $\framme{F}$
      and $\framme{F}'$. Now, suppose that $\phi$ is a formula that is
      satisfied in $\framme{F}$ and $\framme{F}'$ and prove that $\phi$
      belongs to $\logic{L}_2$. It suffices to prove that $\phi$ is valid in
      the canonical model of $\logic{L}_2$.

      First, if $w$ is a world of $\model{M}_{\logic{L}_2}$, then $w \models
      \square \square \psi$ for every formula $\psi$ and thus the
      subframe $R_{\logic{L}_2}^\omega w$ of $\struc{W_{\logic{L}_2}, R_{\logic{L}_2}}$ 
      generated by $w$ is equal to $R_{\logic{L}_2}w$ and $(w,w) \not\in
      R_{\logic{L}_2}$. 

      Then, since $w \models \neg (\diamond \phi \wedge \diamond \neg \phi)$
      for every formula $\phi$, it follows that $\card{R_{\logic{L}_2}}\leq
      1$. Otherwise there are two worlds $t$ and $s$ in $R_{\logic{L}_2}w$. Since
      we work in the canonical model, it means that there is a formula $\psi$
      such that $t(\psi)\neq s(\psi)$. 
      We deduce that $w\not\models\neg(\diamond \psi \wedge \diamond \neg
      \psi)$.

      Now, with the help of the axiom $\square (p \vee \neg p)$, we obtain
      that $R^\omega_{\logic{L}_2}w$ is an $n+1$-$\pi$-morphic image of
      $\framme{F}$ or $\framme{F}'$ and thus $\model{M}_{\logic{L}_2}, w
      \models \phi$.

      We leave to the reader the task to prove, similarly as in (1), that  the
      logic $\logic{L}_2$ is not \textsc{Kripke} complete.

      For (3) and (4), it suffices to mimic the classical proofs.
 \end{proof}
The preceding proposition illustrates the difference between frame definability and
$n+1$-frame definability and gives a first simple example of correspondence
theory for modal many-valued formulas and $n+1$-frames. We should study in more
details  in a forthcoming paper   this dissimilarity in frame definability.


\subsection{Infinitary modal many-valued systems}
Another way of obtaining completeness results is to extend the modal systems
that we have defined by an infinitary deduction rule.
\begin{defn}
  The \emph{infinitary modal many-valued system} has the set $\Form$ of
  well-formed formulas, the \emph{modus ponens}, the necessitation rule and
  the rule
\[
\regle{Inf} \quad \begin{array}{c}
 \phi\oplus \phi,\, \phi \oplus \phi^2,\, \ldots\, ,\,  \phi \oplus \phi^n, \,
 \ldots \\ \hline
 \phi
\end{array}
\]
as deduction rules.

 If $\Gamma$ is a set of axioms and $\phi$ is a formula, we write
$\Gamma \proves_\infty \phi$ if $\phi$  appears in a possibly infinite sequence
$(\psi_\beta)_{\beta \leq \alpha}$ of formulas that belongs to $\Gamma$, are
obtained by substitution in a formula belonging to $\Gamma$ or that are
obtained by the application of a deduction rule from previous formulas of the sequence.
\end{defn}


Note that if $\phi$ is a formula and $\Gamma \proves_{\logic{K}} \phi$ then
$\Gamma \proves_\infty \phi$. We can then state the following completeness result.
\begin{prop}
If $\Gamma \cup \{\phi\}$ is a set of formulas, then $\Gamma
\proves_{\infty} \phi$ if and only if $\phi$ is true in every
many-valued \textsc{Kripke} model of $\Gamma$.    
\end{prop}
\begin{proof}
  If $\Gamma \proves_{\infty}\phi$, the result follows from that fact
  that there is no infinitely great element in the MV-algebra $[0,1]$, so that
  the rule $\regle{Inf}$ preserves tautologies.

  Now, if $\phi$ is true in any model of $\Gamma$, then $\phi$ is infinitely great in the
  algebra $\free_\logic{K}/\struc{\Gamma}$, which means that for any integer
  $m \geq 2$, the element $\phi \oplus \phi^m$ is equal to $\struc{\Gamma}$ in
  $\free_\logic{K}/\struc{\Gamma}$, or equivalently  that $\Gamma \proves_{\logic{K}} \phi \oplus \phi^m$  for any
  integer $m \geq 2$. We can conclude using the rule $\regle{Inf}$.
\end{proof}



\section{Conclusions}
We propose some tracks than could be followed in the future.
\begin{itemize}
  \item Infinitary vs finitary modal systems. The general completeness result
  that is proposed in this paper involves an infinitary deduction rule. Since
  we can get rid off this rule in the case of the finitely-valued modal
  logics, the question to determine the minimal extensions of $\logic{K}$ for
  which the completeness result~--~without the infinitary rule~--~can be stated
  should be considered. We do not know if $\proves_{\infty} \phi$ is
  equivalent to $\proves_{\logic{K}} \phi$.
  \item Varieties of MMV-algebras and MMV${}_n$-algebras. We have not given
  any significant information about these varieties. A good tool for the
  studies of varieties of  MMV${}_n$-algebras could be the topological duality
  constructed in \cite{BruDual}. A problem that could be solved in this way is
  the characterization of finitely generated algebras  
  (following some ideas of \cite{Grigo}).
  \item \textsc{Kripke}-completeness, $n+1$-\textsc{Kripke}-completeness and
  correspondence theory: we should study in detail, with the tool of universal
  algebra and canonical extensions, the dissimilarity between
  \textsc{Kripke}-completeness and $n+1$-\textsc{Kripke}-completeness and
  consider the problem of the correspondence between modal many-valued
  formulas and first order sentences on frames and $n+1$-frames.  
  \item Temporal logic and propositional dynamic logic: we should give interest to
  the construction of some particular systems of modal $n+1$-valued logics such
  as $n+1$-valued temporal logic or dynamic logic.  By their nature, such
  systems could be useful for computer scientists for example.
  \item Extension to more general languages. The results of this paper are
  obtained for the basic modal language with one unary modal operator. They
  should be extended to languages containing $k$-ary modalities.
\end{itemize}
\bibliographystyle{abbrv}


\begin{thebibliography}{10}

\bibitem{Vene1}
P.~Blackburn, M.~de~Rijke, and Y.~Venema.
\newblock {\em Modal logic}, volume~53 of {\em Cambridge Tracts in Theoretical
  Computer Science}.
\newblock Cambridge University Press, Cambridge, 2001.

\bibitem{Chagrov}
A.~Chagrov and M.~Zakharyaschev.
\newblock {\em Modal logic}, volume~35 of {\em Oxford Logic Guides}.
\newblock The Clarendon Press Oxford University Press, New York, 1997.
\newblock Oxford Science Publications.

\bibitem{Chang1}
C.~C. Chang.
\newblock Algebraic analysis of many valued logics.
\newblock {\em Trans. Amer. Math. Soc.}, 88:467--490, 1958.

\bibitem{Chang2}
C.~C. Chang.
\newblock A new proof of the completeness of the {\l} ukasiewicz axioms.
\newblock {\em Trans. Amer. Math. Soc.}, 93:74--80, 1959.

\bibitem{Cign}
R.~L.~O. Cignoli, I.~M.~L. D'Ottaviano, and D.~Mundici.
\newblock {\em Algebraic foundations of many-valued reasoning}, volume~7 of
  {\em Trends in Logic---Studia Logica Library}.
\newblock Kluwer Academic Publishers, Dordrecht, 2000.

\bibitem{Grigo}
L.~Esakia and R.~Grigolia.
\newblock The criterion of {B}rouwerian and closure algebras to be finitely
  generated.
\newblock {\em Polish Acad. Sci. Inst. Philos. Sociol. Bull. Sect. Logic},
  6(2):46--52, 1977.

\bibitem{Fit2}
M.~Fitting.
\newblock Many-valued modal logics. {II}.
\newblock {\em Fund. Inform.}, 17(1-2):55--73, 1992.

\bibitem{Fit1}
M.~C. Fitting.
\newblock Many-valued modal logics.
\newblock {\em Fund. Inform.}, 15(3-4):235--254, 1991.

\bibitem{MvAlgUniv}
J.~Gispert and D.~Mundici.
\newblock M{V}-algebras: a variety for magnitudes with {A}rchimedean units.
\newblock {\em Algebra Universalis}, 53(1):7--43, 2005.

\bibitem{Got}
S.~Gottwald.
\newblock {\em A treatise on many-valued logics}, volume~9 of {\em Studies in
  Logic and Computation}.
\newblock Research Studies Press Ltd., Baldock, 2001.

\bibitem{Tar1}
B.~J{\'o}nsson and A.~Tarski.
\newblock Boolean algebras with operators. {I}.
\newblock {\em Amer. J. Math.}, 73:891--939, 1951.

\bibitem{Tar2}
B.~J{\'o}nsson and A.~Tarski.
\newblock Boolean algebras with operators. {II}.
\newblock {\em Amer. J. Math.}, 74:127--162, 1952.

\bibitem{Kripke}
S.~A. Kripke.
\newblock Semantical analysis of modal logic. {I}. {N}ormal modal propositional
  calculi.
\newblock {\em Z. Math. Logik Grundlagen Math.}, 9:67--96, 1963.

\bibitem{Luka1}
J.~{\L}ukasiewicz.
\newblock O logice tr\' ojwarto'sciowej.
\newblock {\em Ruch Filozoficny}, 5:170--171, 1920.

\bibitem{Luka2}
J.~{\L}ukasiewicz.
\newblock {\em Selected works}.
\newblock North-Holland Publishing Co., Amsterdam, 1970.
\newblock Edited by L. Borkowski, Studies in Logic and the Foundations of
  Mathematics.

\bibitem{Ost}
P.~Ostermann.
\newblock Many-valued modal propositional calculi.
\newblock {\em Z. Math. Logik Grundlag. Math.}, 34(4):343--354, 1988.

\bibitem{BruDual}
B.~Teheux.
\newblock A duality for the algebras of a {\l}ukasiewicz $n+1$-valued modal
  system.
\newblock {\em To appear in Studia Logica}, 2006.

\end{thebibliography}

\end{document}